\documentclass{amsart}
\usepackage{amssymb, euscript}

\begin{document}

\address{Azer Akhmedov, Department of Mathematics,
University of California,
Santa Barbara, CA, 93106, USA}
\email{akhmedov@math.ucsb.edu}

\begin{center} {\bf A NEW METRIC CRITERION FOR NON-AMENABILITY \ {\large I}} \end{center}

\vspace{0.6cm}

\begin{center} {\bf Azer Akhmedov} \end{center}

\vspace{1cm}

\begin{center} {\small ABSTRACT: By studying the so-called traveling salesman groups, we obtain a new metric criterion for non-amenability. As an application, we give a new and very short proof of non-amenability of free Burnside groups with sufficiently big odd exponent.}  \end{center}

\vspace{1.3cm}

    Amenable groups have been introduced by John von Neumann in 1929 in connection with Banach-Tarski Paradox, although earlier, Banach himself had understood that, for example, the group $\mathbb{Z}$ is amenable. J.von Neumann's original definition states that a discrete group is amenable iff it admits an additive invariant probability measure. In 1950's Folner gave a criterion which, for finitely generated groups, can be expressed in terms of the Cayley metric of the group. Using this criterion it is very easy to see that abelian groups are amenable and non-abelian free groups (e.g. $\mathbb{F}_2$) are not amenable.

    \medskip

   In 1980's Grigorchuk [Gr1] introduced a new metric criterion for amenability in terms of the co-growth as a refined version of Kesten's criterion [K]. Using this criterion, A.Ol'shanskii [Ol1] constructed a counterexample to von Neumann Conjecture, and using the same criterion S.Adian [Ad1] proved that free Burnside groups of sufficiently large odd exponent are non-amenable. Recently, A.Ol'shanskii and M.Sapir [OS] found a finitely presented counterexample to von Neumann Conjecture where they establish non-amenability using the same criterion again.

      \medskip

    Based on our studies of traveling salesman groups, we introduce a different metric criterion (sufficient condition) for non-amenability. Roughly speaking, this criterion is based on verifying that certain type of words in a group do not have small length. This, in a lot of cases, seems manageable to show,  compared to measuring the boundary of a finite set. More than that, it was proved(observed) by von Neumann that if a group contains a copy of $\mathbb{F}_2$ then it is non-amenable. But if the  group contains no subgroup isomorphic to  $\mathbb{F}_2$ then it is usually hard to prove non-amenability. However, one can continue in the spirit of von Neumann's observaton. For example, from Folner's criterion it easily follows that if a group is quasi-isometric to a group which contains a copy of $\mathbb{F}_2$ then the group is still non-amenable. Our criteria generalizes much further in this direction.

    \medskip

    {\bf Traveling salesman $(TS)$ groups} were introduced in [Ak1] and [Ak2]. These groups are very useful and play an important role in the constructions showing that the property of containing free subgroup, containing free subsemigroup, satisfying a law, or a girth type, and numerous other properties are not invariant under quasi-isometry. (see [Ak3]).

    \medskip

    Negatively curved groups turn out to be $TS$ (see below). Here, the term "negatively curved" is used in a more general sense than being non-elementary word hyperbolic. On the other hand, $TS$ implies non-amenability, although the class of $TS$ groups form a proper subclass of all non-amenable groups.

     \bigskip

    {\bf Definition of Amenability.} There exists at least 10 wellknown definitions of the notion of amenable group which look quite different from each other but, amazingly, turn out to be all equivalent. There are also related notions like weak amenability, inner amenability, uniform non-amenability, amenable algebra, amenable action, etc. The equivalences between different definitions of amenable group are often respectable theorems.

    \medskip

     The only definition we will be using in this paper is the one obtained from Folner's Criterion. This definition, again, in its own turn, has at least 10 different versions which are all equivalent (and the equivalences in this case are easy exercises).

     \medskip

    {\em Definition 0.} Let $\Gamma $ be a finitely generated group. $\Gamma $ is called amenable if for every $\epsilon > 0$ and finite subsets $K, S\subseteq \Gamma $, there exists a finite subset $F\subseteq \Gamma $ such that $S\subseteq F$ and $\frac {|FK\backslash F|}{|F|} < \epsilon $.

   \medskip

   The set $F$ is called $(\epsilon , K)$-Folner set. Very often one uses the loose term "Folner set", and very often one assumes $K$ is fixed to be the symmetrized generating set. The set $\{x\in F \ | \  xK^{-1}K\subseteq F\}$ will be called {\em the interior of $F$} and will be denoted by $Int_KF$. The set $F\backslash Int_KF$ is called {\em the boundary of $F$} and will be denoted by $\partial _KF$.

   \bigskip

   Despite it's great theoretical value, Folner's Criterion is often very unpractical. For example, to establish non-amenability, in the direct application of the criterion, one has to argue somehow that given a small $\epsilon > 0$, no set can be an $\epsilon $-Folner set, and it is often very hard to rule out all sets at once since Folner sets may have very tricky, versatile and complicated geometry and combinatorics. Another major problem is that, given $\epsilon > 0$, the definition does not demand any upper bound on the size of the minimal $\epsilon $-Folner set. Yet another issue is that the size $|F|$ of $\epsilon $-Folner set is not any explicit function of the size $|\partial F|$.

 \bigskip
 
 {\em Acknowledgement:} I would like to thank to S.Adian, V.Guba, S.Ivanov, A.Olshanski, M.Sapir and W.Thurston for the useful discussions related to the content of this paper. Non-amenability of $TS$ groups was pointed out to me by W.Thurston.

    \vspace{1cm}

       \begin{center}  STRUCTURE OF THE PAPER \end{center}

       \bigskip

       In Section 1, we mention definitions and some facts about $TS$ groups. We also include the proposition saying that $TS$ implies non-amenability (due to W.Thurston).

       In Section 2, we introduce property $(P)$ and prove that it implies $(TS)$. Then we weaken property $(P)$, replacing it with $(P')$ and prove that $(P')$ still implies $(TS)$. Then we formulate a general (i.e. even weaker) property \ $(T_{p,q})$ \ which still implies $(TS)$.

       In Section 3, we discuss introduce property  $(TS)(\lambda )$ and discuss it for products of groups.

       Finally, in Section 4, we discuss some applications of the criteria that we introduce in earlier sections.

      \medskip

      The proofs are in somewhat hierarchic style, that is several times we prove a result stronger than the previous result. We do this to make the proofs more readable.

  \vspace{1cm}

  \begin{center} SOME NOTATIONS.  \end{center}

  \bigskip

 {\bf 1.} The Cayley metrics in all finitely generated groups will be assumed to be left invariant, and denoted as $|.|$ or $d(.,.)$; it will be always clear from the context which group we are considering. The balls of radius $r>0$ centered at the element $g$ will be denoted by $B_r(g)$.

\medskip

 {\bf 2.} By ternary trees (with origin)  we mean a finite tree of which all vertices have either valence \ $1$ - end vertex or $3$ - internal vertex (or middle vertex). One of the internal vertices (or end vertex, if the tree has only one vertex) is marked as an origin. For each vertex \ $v$ \ in a ternary tree, the edge which joins this vertex with the vertex closer to the origin will be called {\em upper edge}. The other two edges will be called the {\em lower edges}. By easy induction one sees that the number of end vertices of a ternary tree is two more than the number of internal vertices. We call a finite (infinite) tree a $ray$ if all vertices have valence 2, except two (one) which we will call the ends of a ray. The vertices of infinite rays will sometimes be denoted by integers $0, 1, 2, 3, \ldots $

 \medskip

 {\bf 3.} By a $plane \ ternary \ tree \ with \ origin$, we mean a ternary tree with origin embedded in a plane with a cartesian coordinate system such that the origin is at the point $(0,0)$ and if the vertex $x$ is at the distance $n$ from the origin then $x$ lies on the horizontal line $x = -n$. We will put an order on the set of vertices of this graph in the following manner: the vertex $x$ is bigger than line the vertex $y$ if the distance from $x$ to the origin is bigger than the distance from $y$ to the origin, and if $x$ and $y$ are on the same distance from the origin, then $x>y$ if $x$ is on the left side of $y$. The vertex at the distance $n$ from the origin will be called a vertex of level $n$.

\medskip

 {\bf 4.} A sequence is called aperiodic if it has no two consecutive (adjacent) identical subsequences. It is called $k$-aperiodic if there exist no \ $k+1$\  consecutive identical subsequences. For example, the sequence \ 1213123213 \ is aperiodic, and the sequence \ 1213132131 \ is 2-aperiodic. An element $\xi $ in the free group $ \mathbb {F}_k, (k\geq 1)$ with a generating set $\{a_1,  \ldots , a_k\}$ is called aperiodic if it is a reduced word and aperiodic in the alphabet $\{a_1, a_1^{-1}, \ldots  , a_k, a_k^{-1}\}$.

 \medskip

 {\bf 5.} We will denote a cardinality of a set $A$ by $|A|$.

 \medskip

  {\bf 6.} For a finitely generated group $G$ with a fixed Cayley metric $d$, a closed path $\beta $ is a function $\beta : \{0, 1, 2, \ldots n\} \rightarrow G$ where $d(\beta (i), \beta (i+1)) = 1$ for all $0\leq i \leq n-1$, $\beta (0) = \beta (n)$. $n$ will be called the length of $\beta $ and be denoted as $n = l(\beta )$. For any set $A\subseteq G$ we will denote $N(\beta ; A) = |\{ i\in \{0, 1, \ldots n\} \ | \ \beta (i)\in A\}|$.

\medskip

{\bf 7.} In the free group $ \mathbb {F}_k, (k\geq 1)$ with a standard generating set $S = \{a_1, a_1^{-1} \ldots , a_k, a_k^{-1}\}$ and left invariant Cayley metric with respect to this generating set $S$, we say a reduced word $E$ is a subword of a (not necessarily reduced) word $W$ if there exists reduced words $X, Y$ such that $W = XEY$ and $|W| = |X| + |E| + |Y|$.  In this case we also write $W\equiv XEY$ and say $E$ occurs in $W$. Any of the words $X, E, Y$ is allowed to be empty. If $W\equiv XA^kY$ then $A$ is called a periodic subword of order at least $k$. If $A = c_1c_2\ldots c_n, \ c_i\in S, 1\leq i\leq n$ and
$W = d_1\ldots d_pc_1c_2\ldots c_nf_1\ldots f_q, \ d_i, f_j\in S, 1\leq i\leq p, 1\leq j\leq q$ are reduced words, and $m < q$ is a positive integer, then $R_m(A;W)$ denotes the subword $c_{1+m}\ldots c_nf_1\ldots  f_m$. So, $R_m(A;W)$ is the shift of $A$ by $m$ units to the right.

\vspace{1cm}

      \begin{bf} 1. Definitions and basic properties. \end{bf}

  \bigskip

         We borrow some definitions from [Ak2]

For a (usually locally finite) graph $X$ and its finite subset of vertices $Y$, the traveling salesman problem usually deals with the minimal path in $X$ which passes through all vertices in $Y$. The length of such a closed path will be denoted by $L(Y)$. One of the usual problems is to measure the ratio \ $|Y|/L(Y)$.

\medskip

 {\em Definition 1}. Let \ $A$ \ be any set in a group \ $\Gamma $. We call a set $B$ in $\Gamma  \ A$-related if for any \ $x\in B$, \ there exists \ $y\in B$ such that \ $|\{y^{-1}x, x^{-1}y\}\cap A| \geq 1$ \ . $y$ will be called an  $A$-neighbour of $x$.

\medskip

 {\em Definition 2}. A group \ $\Gamma $ \ is said to belong to the class \ $TS$ \ if for any \ $\lambda > 1$, \ there exists \ $\xi \in \Gamma $ \ such that for any \ $\{\xi \}$-related finite set \ $S$ in $\Gamma , \ L(S) \geq \lambda |S|$.

\bigskip

 {\bf Independence on the choice of generating set:} In the above definitions, and everywhere in this paper, we assume that \ $\Gamma $ \ has a fixed generating set which defines the length in \ $\Gamma $ by left invariant Cayley metric. But since two metrics corresponding for different choice of generating sets are bi-Lipshitz equivalent, it is clear that the property $TS$ is independent on this choice.

\medskip

 {\bf Revising the set S:} {\em Definition 3.} Let $\Gamma $ be a finitely generated group, $\xi \in \Gamma $, and let $S$ be a $\xi $-related finite subset of $\Gamma $. A subset $S'\subseteq S$ is called a revision of $S$ with respect to $\xi $, if the following conditions are satisfied:

  (i) $|S'| \geq \frac {2}{3}|S|$.

  (ii) there exists a subset $S_0\subset S$ such that $S' = \sqcup _{x\in S_0}\{x, x\xi \}$. (note that $\sqcup $ denotes a disjoint union). If $y\in \{x, x\xi \}$ for some $x\in S_0$, then the only element of the set $\{x, x\xi \}\backslash \{y\}$ is called a $\xi $-neighbor of $y$ and will be denoted by $N_{\xi ,S}(x)$. If the set $S$ and the element $\xi $ are fixed, then we will sometimes drop it in the index by writing $N(x)$.

\medskip

 {\bf Remark 1.} Clearly, every finite $\xi $-related set $S$ has a revision. In other words, $S$ can be revised. If $\xi $ is fixed, we will often drop the phrase {\em with respect to $\xi $}.

 \medskip

 {\em Definition 4.} If $S$ is a revision of $S$ (i.e. of itself) then we will say $S$ is revised.

 \medskip

 One of the major goals is to understand the class of $TS$ groups. For example, it is very easy to see that the group $\mathbb {Z}^n$ is not $TS$, because an $n$-dimensional parallelogram can be traveled in at most twice the cardinality of the parallelogram. The following propositions give more insight.

  \bigskip

    {\bf Proposition 1.} Assume \ $\Gamma $ \ is a finitely generated group which has a non-trivial central
    element of infinite order. Then \ $\Gamma $ \ is not $TS$.

     PROOF: See [Ak2].

     \medskip

     {\bf Corollary 1.} Finitely generated nilpotent groups are not $TS$ groups.

     \bigskip

     {\bf Proposition 2.} \ $\mathbb {F}_k$ \  is  a  $TS$ group for any \ $k\geq 2$.

      PROOF: See [Ak2]

      \bigskip

      {\bf Remark 2.} For a finite set $S$ in a finitely generated group $G$ with a fixed Cayley metric $d$,
      we define $L'(S) = inf _{\beta } (l(\beta ) - N(\beta ; S))$ where infimum is taken over all closed paths $\beta $ which passes through every element of $S$. Then, $L'(S) < L(S)$. From the proof of Proposition 2, it is obvious that we indeed prove slightly more: for any aperiodic word $\xi $ of length $[4\lambda + 1]$ and for any $\xi $-related finite set $S$, the inequality $L'(S) > \lambda |S|$ holds.

      \medskip

      The general rough intuition is that negatively curved groups are $TS$ groups, and groups which are far from being negatively curved are not.

      \medskip

      {\bf Proposition 3.}(W.Thurston) \ Finitely generated amenable groups are not $TS$.

   \medskip

       PROOF:  Assume \ $\Gamma $ is $TS$ and amenable. For \ $\lambda = 3$, \ let \ $\xi \in \Gamma $ \ be such that for any $\xi $-related set \ $S$, \ $L(S) > 3|S|$.  Let \ $F$ \ be a \ $(B_{|\xi |}(1),
\frac {1}{100})$-invariant connected Folner set, and denote \ $\partial F = \{x\in F : \{x\xi ,x\xi ^{-1}\}\cap F = \emptyset \}$. Then \ $F \backslash \partial F$ \ is \ $\xi $-related. On the other hand, \ $F$ \ is a finite connected subgraph of a Cayley graph of \ $\Gamma $, so it has a spanning tree $T$ (so $|T| = |F|$). Then \ $L(F\backslash \partial F) \leq L(F) \leq 2|T| = 2|F| \leq 2.5|F\backslash \partial F|$. Contradiction. $\square $

    \bigskip

     {\bf Remark 3.} The fact about the relation of spanning tree to traveling salesman problem was observed even in 1950's by J.Kruskal (see [K]).

     \medskip

     {\bf Remark 4.} It is easy to see that the direct product of an infinite non-$TS$ group with any group is not a $TS$ group. So, for example, \ $\mathbb{Z}\times \mathbb{F}_2 \notin TS$. This shows that the class of $TS$ groups is a proper subclass of the class of all non-amenable groups.

     \bigskip

     \vspace{1cm}

     {\bf 2. Property ($P$) and Property ($P'$)}

     \bigskip

     The following propositions justify our intuition about negatively curved groups.

     \medskip

 {\bf Proposition 4}. Every torsion free word hyperbolic group which is not virtually cyclic is a $TS$ group.

 \medskip

 {\bf Remark 5:} It is possible to drop the condition of being torsion free but we do not discuss here in this generality. In connection with this, we would like to remind that it is not known if every word hyperbolic group contains a finite index torsion free subgroup.

\medskip

 Although Proposition 4. is much more general than Proposition 2, still the proof of Proposition 2. is interesting and the idea is useful in some other situations (e.g. see Remark 11.)

\bigskip

 In order to prove  Proposition 4.  we will consider the following property for finitely generated groups.

\medskip

 $(P)$: A finitely generated nontrivial group $K$ is said to satisfy the property $(P)$ if for any \ $r > 0$, there exists \ $\xi \in K$ \ such that $$|\xi ^{\epsilon _1 }x_1\xi ^{\epsilon _2 }x_2\ldots \xi ^{\epsilon _k }x_k| > r \ (**)$$ for any positive $k$ and \ $x_1,\ldots ,x_k \in K, \ x_i \neq 1, \ |x_i| \leq r, \ \epsilon _i\in \{+1, -1\}$

\medskip

 We also introduce the following properties:

 $(P')$: A finitely generated nontrivial group \ $K$ \ is said to satisfy a property $(P')$ if for any \ $r > 0$, \ there exists \ $\xi \in K$ \ such that for every aperiodic sequence \ $x_1, x_2, \ldots x_k$ of elements in \ $K$, the inequality \ $|\xi ^{\epsilon _1}x_1\ldots \xi ^{\epsilon _k}x_k| > r$ \ satisfied for every \ $k\in \mathbb{Z}_+, \ x_i \neq 1, \ |x_i| \leq r, \ \epsilon _i \in \{+1, -1\}$.

\medskip

  The next is some modification of property $(P')$:

  \medskip

  $(P_n')$. Let $n$ be a positive integer. A finitely generated nontrivial group \ $K$ \ is said to satisfy a property $(P_n')$ if for any \ $r > 0$, \ there exists \ $\xi \in K$ \ such that for every $n$-aperiodic sequence \ $x_1, x_2, \ldots x_k$ of elements in \ $K$, the inequality \ $|\xi ^{\epsilon _1}x_1\ldots \xi ^{\epsilon _k}x_k| > r$ \ satisfied for every \ $k\in \mathbb{Z}_+, \ x_i \neq 1, \ |x_i| \leq r, \ \epsilon _i \in \{+1, -1\}$.   (So $(P_1')$ is equivalent to $(P')$)

\bigskip

 Now we are going to prove that both of the the properties \ $(P)$ and \ $(P')$ \ are sufficient to imply the property $TS$. The proof for the case of $(P)$ is easier, and the case of $(P')$  will be a refinement of the case for $(P)$. Indeed, we prove that $(P_{10}' )\Rightarrow (TS)$, which is a slightly weaker result (property  $(P_{10}' )$ is (slightly) stronger than $(P')$).

\bigskip

 We will need the following

 \medskip

 {\em Definition 5.} Let $S\subseteq \Gamma $ be a finite revised subset with respect to $\xi \in \Gamma $. A finite ternary tree $T$ with origin is called {\em $(S, \xi )$-tree} if

 (i) all internal vertices consists of triple of elements of $S$, all end vertices are either pairs of elements of $S$, or are singletons of $S$.

 (ii) if $u, v$ are adjacent vertices of $T$ of level $n-1$ and $n$ respectively, $n\in \mathbb{N}$, then there exist $x, y\in S$ which belong to vertices $u$ and $v$ respectively such that $y = N_{\xi }(x)$.

 (iii) vertices of $T$ are disjoint as subsets of $S$.

 \bigskip
 
  When it is clear from the context, instead of $(S,\xi )-tree$ we will simply say {\em tree}.
   
\bigskip

 {\bf Case $(P)$:} Let $\lambda > 1,  \ r > 12\lambda $, and \ $S$ be a $\xi $ related finite set where $\xi $ is chosen to satisfy the inequality $(**)$.

\medskip

 Assume  $\alpha $  is a shortest closed path passing through all vertices of $S$, and let $\psi = \{x_1, x_2, \ldots ,x_n \}$ be a sequence of elements of $S$ given in the order along $\alpha $. Clearly $n \geq |S|$. \ We divide the sequence $\psi $ into the \ $l$ \ subsequences \ $\{x_1, \ldots x_{i_1}\}, \{x_{i_1+1},\ldots ,x_{i_2}\}, \\ \ldots , \{x_{i_{l-1}+1}, \ldots x_n\}$ where $d(x_j,x_{j+1}) \geq r/2$ if and only if $j \in \{i_1,i_2,\ldots ,i_{l-1}, n\}$. Now we will delete elements from these subsequences as follows: for every $j\in \{2,3,\ldots n\}$, \ if $x_j$ belongs to one of the subsequences but \ $x_j = x_i$ for some $i < j$, then we delete $x_j$ from this subsequence [{\em so this is some sort of revision but not the revision in the sense of Definition 3.}]. Thus we obtain subsequences(some of them maybe empty) of the previous subsequences. If necessary we divide each of them again into more subsequences and obtain subsequences (which we will call $ pieces $) \ $\{z_1, \ldots z_{j_1}\}, \ldots ,\{z_{j_{l-1}+1}, \ldots ,z_m\}$ \ such that $m = |S|$, \ and \ $d(z_j,z_{j+1}) > r/2$ \ if and only if \ $j \in \{j_1,\ldots ,j_{l-1},m\}$. In addition, since $\alpha $ is the minimal path we also have \ $\sum _{j=1}^{m}d(z_j,z_{j+1}) = L(S)$.

\medskip

 We divide each piece \ $\{z_p, \ldots ,z_q\}$ \ into $ short \ segments $ \ $\{z_p,z_{p+1},z_{p+2}\}, \{z_{p+3},z_{p+4},z_{p+5}\}, \\ \ldots $, where each short segment consists of 3 elements except the last may have 1 or 2 elements. We will call the short segments of length 3 (less than 3) a $ normal $($ incomplete $) segment.

\medskip

 Now let \ $z$* $\in S$. We can and will assume that \ $S$ \ is revised. We will associate to $z$* a finite {\em plane ternary $(S, \xi )$-tree $T_1$ with origin} as follows: each vertex of \ $T_1$ \ will be a short segment, which will be normal for origin or middle vertices, and incomplete for end vertices. The origin will consist of a 3-point set \ $s_0 := \{z$*,$z',z''\}$ \ where \ $z'$ \ and \ $z''$ \ are the two other elements in a short segment containing \ $z$*. (If this segment is incomplete then \ $T_1$ \ will be a trivial one-point binary tree). Let \ $z''', z'''', (z$*$)'$ \ be the \ $\xi $-neighbors of \ $z', z''$, \ and $z$* with  \ $s_1, s_2, s'$ \ being the short segments containing \ $z''', z'''', (z$*$)'$ \ respectively. The three adjacent vertices to \ $s_0$ \ will be exactly \ $s_1, s_2$ and $s'$. If any of the \ $s_1, s_2, s'$ \ are incomplete then that will be an end vertex, otherwise we continue the process until it becomes impossible. We finish this way the construction of \ $T_1$, \ and then taking \ $z$** $\in  S \backslash T_1$ we start a new tree and so on, thus covering \ $S$ \ with plane ternary $(S, \xi )$-trees \ $T_1, T_2, \ldots ,T_m$. The trees are built inductively, i.e. for each $1\leq i\leq m, n\in \mathbb{N}$ we do not start any vertex of level $n+1$ of $T_i$ until we have built all vertices of level $n$ of $T_i$, and we do not start building the tree $T_{j+1}, 1\leq j\leq m-1$, until we have finished building the tree $T_j$.

\medskip

 By construction, the vertex sets (which are short segments) of the same tree are always disjoint, and the vertex sets of two different trees may intersect only if they are both ends, thus an end vertex may belong to at most two trees.
 By their definition and by property $(P)$, an element of any vertex is in a distance at least \ $r$ \ from the elements of any other vertex of the same tree. [because, by construction, if $\zeta _1, \zeta _2$ are elements of $\Gamma $ which belong to the different vertices of the same tree $T_i, 1\leq i\leq m$, then $\zeta _1 = \zeta _2W$ where $W = y'\xi ^{\epsilon _1 }x_1\xi ^{\epsilon _2 }x_2\ldots \xi ^{\epsilon _k }x_k,  \ x_1, \ldots , x_k \in B_{r/2}\backslash \{1\}, y'\in B_{r/2}, \epsilon _i\in \{-1, 1\}, 1\leq i\leq k$, thus $|W| \geq r$.]

 Besides, by construction, every end vertex contains an element in a distance more than \ $r/2$ \ from an element next to it in the order of the path \ $z_1, z_2, \ldots ,z_m $.

\medskip

 Now assume \ $card(T_i) = t_i, \ 1\leq i \leq m$. Then the number of end elements is at least \ $\sum _{1\leq i\leq m}\frac {t_i}{2} = \frac {1}{2}\sum _{1\leq i\leq m} t_i > \frac {1}{6} |S|$. On the other hand, because of the above remarks one has \ $L(S) > r/2 \sum _{1\leq i\leq m}\frac {t_i}{2} > r/12 |S| > \lambda |S|$. This finishes the proof for case $(P)$.

\bigskip

 For the case of \ $(P_{10}')$ \ we will use the following combinatorial lemmas.

\medskip

 {\bf Lemma 1:} Let \ $T$ \ be an infinite ternary tree with origin. Then one can assign one of the three letters $A, B$, and $C$ to each edge such that a sequence along any simple path in $T$ is 3-aperiodic, i.e. one does not have 4 consecutive identical subsequences.

\medskip

 Proof: The lemma is true for a ray, i.e. there exists an aperiodic sequence in three letters $A, B$, and $C$(see [Ad], page 5). Let \ $\omega $ \ be an example of a such sequence. Then we take all infinite rays in $T$ which start from the origin and label these rays identical to \ $\omega $.

 Any simple path \ $L$ \ is a union of at most two simple pathes \ $L_1$ and $L_2$ such that each piece $L_i$ is contained in a ray starting from the origin. Therefore \ $L$ \ is 3-aperiodic.   $\square $

\bigskip

 {\bf Lemma 2:} Let \ $T$ \ be a plane ternary tree with origin, and \ $F$ \ be a finite set, such that for each vertex \ $v\in T$, \ there is a subset \ $F_v \subset F$ \ assigned to \ $v$, \ with \ $|F_v| \geq 4$. \ Then for every \ $v\in T$ \ one can choose \ $x_v, y_v, \in F_v$ \ and label with them the two lower edges adjacent to \ $v$ \ s.t. every simple path in \ $T$ \ is 10-aperiodic. Moreover, one can do the labeling inductively, more precisely, one first labels the edges adjacent to the root of a tree \ $T$ \ and then for any \ $v\in T$, \ if the upper edge of \ $v$ \ is labeled but the lower edges are not labeled, then one can choose \ $z_v\in F_v$ \ (called inadmissible element at the vertex \ $v$) \ s.t. all elements of \ $F_v\backslash z_v$ \ are admissible, i.e. one can label the lower edges adjacent to \ $v$ \ with any two \ $x_v, y_v \in F_v\backslash z_v $.
\medskip

 We first prove this lemma (indeed slight generalization of it) for a ray.

\medskip

 {\bf Lemma 3.} Let \ $R$ \ be an infinite ray and \ $F$ \ be a finite set such that at each vertex \ $v\in R$ \ there is a subset \ $F_v\subset F, \ |F_v| \geq 2$ \ assigned to $v$. Then, starting from the origin, one can choose \ $z_v\in F_v $ and label the lower edge adjacent to $v$ (i.e. the edge joining vertex $v$ with the vertex further from origin) with any (!) $x_v\in F_v\backslash z_v $ and continue inductively so that the resulting sequence \ $\{x_v\}$ is 4-aperiodic.

\medskip

 Proof: We proceed by induction on $|F|$.

\medskip

 For \ $|F| = 2$, if \ $F = \{a,b\}$ \ let \ $\omega _0 $ \ be an aperiodic sequence in symbols \ $ab, ba, aa, bb$ \ and let us denote by \ $\omega ' = \{\omega '_1, \omega '_2, \ldots \} $ \ the same sequence \ $\omega _0 $ \ but in symbols $a$ and $b$. Then \ $\omega '$ \ is a 4-aperiodic sequence. Then for any $n\geq 1$, we choose \ $z_n = \omega _n'$ \ and then the choice of \ $x_n$ \ is forced: if $z_n = a$ then $x_n = b$, and if $z_n = b$ then $x_n = a$. Clearly, the sequence $\{x_n\}_{n=1, \infty }$ is $4$-aperiodic.

\medskip

 Assume now the lemma is true for  $|F| = k$. For the case  $|F| = k+1$, let \ $F = \{a_1, \ldots a_k, a_{k+1}\}$, \ and assume \ $\{c_1, c_2, \ldots \}$ \ be a sequence of inadmissible elements, and $\theta _k$ be a sequence of 4-aperiodic $k$ symbols which we would choose as the sequences $ \{z_1, z_2, \ldots \}, \ \{x_1, x_2,\ldots \} $ respectively if $x_n \neq a_{k+1}$ for any $n\geq 1$; the existence of such sequences follows from inductive hypothesis.

 At each vertex $n$ of the ray $R$ we choose $z_n = c_n$ until the first time $x_{n-1} = a_{k+1}$. Then

  1) if $\omega '_1 = a$ we set $z_n = a_{k+1}$

  2) if $\omega '_1 = b$ we set $z_n = c_n$ and if also then $x_n = a_{k+1}$ we set $z_{n+1} = a_{k+1}$.

 Then we continue assigning the elements of the sequence \ $\{c_1, c_2, \ldots \}$ to the edges of the ray which was stopped at $c_{n-1}$ (case 1) or $c_n$ (case 2) until again \ $x_{m-1} = a_{k+1}$ for some $m > n$ and we do the same (now things will depend if $\omega '_2 = a$ or $b$) and so on. Since $\omega '$ and $\theta _k$ are both 4-aperiodic the resulting sequence $\{x_n\}$ is 4-aperiodic. Lemma is proved.

\bigskip

 Proof of Lemma 2. Let \ $R_1, R_2, \ldots $ be a sequence of rays starting from the origin such that their union is the whole tree $T$.

  First, we label the ray \ $R_1$ \ and all edges adjacent to it inductively (as in the statement of the Lemma 3). Then any simple path in the labeled subtree is $(4+1) = 5$-aperiodic. Assume we have labeled the rays \ $R_1, \ldots , R_n$ together with all edges adjacent to them. Let $v\in R_{n+1}\backslash \cup _1^n R_i $ be a vertex in $R_{n+1}$ such that all vertices in the segment from the origin up to $v$ (excluding $v$) belong to $\cup _1^n R_i$. Then there exists a ray \ $R_i, \ 1\leq i \leq n$ such that \ $R_i$ and $R_{n+1}$ coincide from the origin up to $v$ (assume this segment has length $d$), excluding $v$ again, and we choose the sequence \ $\{z_d, z_{d+1}, \ldots \}$ of inadmissible elements for the ray \ $R_{n+1}$ \ as in the case of \ $R_i$.

 This way we label all edges of the tree $T$ and it is clear that any simple path will be $5\times 2 = 10$-aperiodic. Lemma is proved.

\bigskip

 {\bf Case $(P_{10}')$:} Let \ $\lambda > 1, r > 96\lambda $, and \ $\xi \in K$, such that the inequality $|\xi ^{\epsilon _1}x_1\ldots \xi ^{\epsilon _k}x_k| > r$ is satisfied for every \ $k\in \mathbb{Z}_+, \ x_i \neq 1, \ |x_i| \leq r, \ \epsilon _i \in \{+1, -1\}$. We are going to show that \ $K$ then satisfies \ $TS(\lambda )$. \ Assume \ $S$ \ is any finite $\xi $-related set in $K$. As in the previous case we divide the set \ $S$ \ into \ $l$ \ pieces \ $\{z_1, \ldots ,z_{j_1}\},\{z_{j_1+1},\ldots ,z_{j_2}\},\ldots ,\{z_{j_{l-1}+1}, \ldots ,z_{j_l}\}$ \ where \ $j_l = |S|$, and \ $d(z_j,z_{j+1}) > r/4$ if and only if \ $j \in \{j_1, j_2, \ldots ,j_l\}$ (here we assume that \ $z_{j_l+1} = z_1$).

\medskip

 We will again divide the set $S$ into with plane ternary $(S, \xi )$-trees $T_1, T_2, \ldots , T_m$ with origin. Again, the trees will be built inductively, i.e. for each $1\leq i\leq m, n\in \mathbb{N}$ we do not start any vertex of level $n+1$ of $T_i$ until we have built all vertices of level $n$ of $T_i$, and we do not start building the tree $T_{j+1}, 1\leq j\leq n-1$, until we have finished building the tree $T_j$. The major differences from the previous case are that we will not do this by the help of dividing pieces into strict short segments, plus, the trees will be disjoint  but an end vertex of one tree may lie in the piece far from the ends of this piece.

\medskip

 Assume  $\alpha $ is a shortest closed path passing through all vertices of $S$, and let $\psi = \{z_1, z_2, \ldots ,z_n \}$ be a sequence of elements of $S$ given in the order along $\alpha $. Again, we assume that the set $S$ is revised. Let \ $z$* $\in S$, and assume \ $z$* belongs to the piece \ $\{z_p,\ldots ,z_q\}$. The origin of the tree \ $T_1$ \ will be a 3-point set \ $\{z^*, z', z''\}$ \ where \ $z', z''$ \ are the two immediate consecutive left neighbors of $z$* in the piece \ $\{z_p, \ldots ,z_q\}$. (If there are no left neighbors of $z$* then we choose $z',z''$ to be the immediate consecutive right neighbors, if there is only one left neighbor then we choose \ $z'$ \ to be this, and \ $z''$ \ to be the immediate right neighbor, and finally, if the piece contains less than three elements then this piece will be the only vertex of $T_1$, and by taking \ $z$** $\in S \backslash T_1$ \ we start building the tree $T_2$.)

 \medskip

 Assume \ $z''', z''''$ \ are the \ $\xi $-neighbors of \ $z'$ and $z''$ respectively. Let \ $\{u_1, u_2, u_3, u_4\}$, \ $\{v_1, v_2, v_3, v_4\}$ \ be the immediate left and right neighbors of $z'''$ in the piece it belongs to. [More precisely, we need to denote the sets by $\{u_1, \ldots u_i\}$, \ $\{v_1, \ldots, v_j\}$ with $1\leq i, j\leq 4$ since we may have less than 4 elements to the right or to the left of $z'''$ in the piece to which $z'''$ belongs].

 \medskip

 By Lemma 2, only one of these at most 8 elements can be inadmissible. [Here, if $v$ denotes the vertex to which $z'''$ belongs, then, for the application of the lemma, one can take the set $F$ to be equal to $S^{-1}S$ and the set $F_v$ to be equal to $\{(z''')^{-1}u_1, (z''')^{-1}u_2, (z''')^{-1}u_3, (z''')^{-1}u_4, (z''')^{-1}v_1, (z''')^{-1}v_2, (z''')^{-1}v_3, (z''')^{-1}v_4\}$.]
 If any of the $v_1, v_2, v_3, v_4$ belong to another vertex of $T_1$ or to the vertex of previously built trees (one may notice that since there are no previously built tree at this point, and we just started to build the tree $T_1$, this condition is void {\em at the moment}), then we will include the nearest(in the order of a piece) two of the remaining admissible elements to the vertex of $z'''$, i.e. the vertex to which $z'''$ belongs. If the number of such vertices is less than two (in particular, if there are no such elements) then we will add the immediate admissible elements among  $u_1, u_2, u_3, u_4$ to this vertex and thus complete the vertex, i.e. to the vertex of $z'''$. We do the same for  $z''''$ and continue the process until it becomes impossible, i.e. we reach all end vertices of the tree  $T_1$.

 \medskip

 More generally, if we are building the vertex $\bar v$ of level $n\geq 1$ of the tree $T_i, 1\leq i\leq m$, where an element $\bar {z}$ in this vertex is a $\xi $-neighbor of some element $\bar {z'}$ from the vertex of level $n-1$, then we let $\{\bar {u_1}, \bar {u_2}, \bar {u_3}, \bar {u_4}\}$, \ $\{\bar {v_1}, \bar {v_2}, \bar {v_3}, \bar {v_4}\}$ \ be the immediate left and right neighbors of $\bar z$ in the piece it belongs to. By Lemma 2, only one of these at most 8 elements can be inadmissible. If any of the $\bar {v_1}, \bar {v_2}, \bar {v_3}, \bar {v_4}$ belong to another vertex of $T_j, j < i$ or to the previously built vertex of $T_i$, then we will include the nearest (in the order of a piece) two of the remaining admissible elements to the vertex of $\bar z$, i.e. the vertex to which $\bar z$ belongs. If the number of such vertices is less than two [in particular, if there are no such elements; by {\em "such vertices"}, we mean admissible elements among $\bar {v_1}, \bar {v_2}, \bar {v_3}, \bar {v_4}$ which do not belong to another vertex of $T_j, j < i$ or to the previously built vertex of $T_i$] then we will add the immediate admissible elements among $\bar {u_1}, \bar {u_2}, \bar {u_3}, \bar {u_4}$ to this vertex and thus complete the vertex.

 \medskip

 Along any path connecting one vertex of the tree $T_1$ to another vertex of $T_1$ one would read a word of the type $W = y\xi ^{\epsilon _1 }x_1\xi ^{\epsilon _2 }x_2\ldots \xi ^{\epsilon _k }x_k$ where $x_1, \ldots , x_k$ is a 10-aperiodic sequence of elements from $B_{r/2}\backslash \{1\}$, $y\in B_{r/2}$, $\epsilon _i\in \{-1, 1\}, 1\leq i\leq k$, and since, by property $(P_{10}')$, these words have length at least $r$, we obtain that the vertices of $T_1$ are indeed disjoint. Then we take $z$** $\in S\backslash T_1$ and start building the tree $T_2$ and so on, thus partitioning the set \ $S$ \ into the trees \ $T_1, \ldots T_m$.

\medskip

By Lemma 2, only one of these at most 8 elements can be inadmissible. If any of the \ $v_1, v_2, v_3, v_4$ belong to another vertex (of $T_1$), then we will include the nearest(in the order of a piece) two of the remaining admissible elements to the vertex of $z'''$, i.e. the vertex to which $z'''$ belongs. If the number of such is less than two (in particular, if there are no such elements) then we will add the immediate admissible elements among \ $u_1, u_2, u_3, u_4$ to this vertex and thus complete the vertex, i.e. to the vertex of $z'''$. We do the same for \ $z''''$ and continue the process until it becomes impossible, i.e. we reach all end vertices of the tree  $T_1$.

 By our construction the vertices of each tree are 3-point set if it is a middle vertex, and 1-point or 2-point set if it is an end vertex. In contrast from the previous case, the elements of end vertices are not always the first or the last 3 elements of pieces but may lie inside the piece far (i.e. in a distance at least 4) from its ends. Let \ $V''$ \ be the set of such elements, and \ $V'$- the set of remaining elements of end vertices. By construction, \ $|V''| \leq 1/8 |S|$. \ On the other hand, \ $|V'| + |V''| > \frac {1}{3} (\frac {1}{2}|S|) = \frac {1}{6}|S|$. \ Then \ $|V'| > \frac {1}{24}|S|$.

\medskip

 Each element of \ $V'$ \ is in at least \ $r/4$ \ distance from the previous or next element to itself (or to the elements following the previous or next element) in the order \ $\{z_1, \ldots ,z_m\}$ \ of the set \ $S$. Therefore \ $L(S) \geq \frac {r}{4}|V'| > \frac {r}{96}|S| > \lambda |S|$. This finishes the proof for case $(P_{10}')$. By slightly more complicated arguments one could also establish $(P')\Rightarrow (TS)$.

\bigskip

 Thus we proved the following

 \medskip

  {\bf Proposition 5.} The property $(P_{10}')$ [and property $(P)$] implies $TS$.

  \bigskip

 {\bf Proof of Proposition 4:} What remains is to show that the torsion free non-virtually cyclic hyperbolic group \ $\Gamma $ satisfies property $(P)$. This is proved in [Ak2].

\bigskip

 {\bf Remark 6.} Notice that property \ $(P)$ \ implies that the group contains free subgroup of rank $2$(the elements of type \ $x\xi x$ and $y\xi y$ will generate free subgroup for any \ $x,y\in B_{r/2}\backslash \{1\}, \ x\neq y, y^{-1}$), However, property $(P')$ does not imply the existence of free subgroup since free Burnside groups satisfy $(P')$. It is interesting therefore to test the property $(P')$ independently for some known examples of amenable groups (the answer should be negative by Prop.3 and Prop.6.) Indeed, property $(P')$ easily implies exponential growth, so the famous Grigorchuk's group [Gr2] does not satisfy it. For solvable groups, by choosing $x_i$'s from the last subgroup in the derived series (which is abelian), we see that they do not satisfy property $(P')$ either. It is also not difficult to show this for all elementary amenable groups.

 \bigskip

  Let us introduce another property for finitely generated groups.

  \medskip

  {\bf Property ($T_{p,q}$):} Let \ $\Gamma $ \ be a finitely generated nontrivial group,\ $\xi \in \Gamma $, \ $T$ \ be an ordered ternary tree with origin, \ $p\geq 4q$. Let one starts a process of assigning each edge an element of the set \ $B_r(1)\backslash \{1\}$ \ and each vertex the element \ $\xi $ or $\xi ^{-1}$ in the following inductive way: for any vertex \ $v_n$ \ we choose $\xi $ or $\xi ^{-1}$ and assign it to $v_n$ and we choose arbitrary set $F_{v_n}$ \ of \ $p$ \ elements in $B_r\backslash 1$, and then we choose some subset \ $U_{v_n}$ \, and then choose some two elements $x_{v_{n}}, y_{v_{n}}$ of \ $F_{v_n}\backslash U_{v_n}$ \ and assign it to the two lower edges of $v_n$. Then along any simple path in the tree one reads a word like \ $\xi ^{\pm }x_1\xi ^{\pm }x_2\ldots \xi^{\pm }x_k$ where $k$ is the length of the path. We say that \ $\Gamma $ \ satisfies property \ $T_{\{p,q\}}$  if for any $r>0$, there exists \ $\xi \in \Gamma $ such that for any $n\in \mathbb {N}$, \ and for any choice of $\xi , \xi ^{-1}$ at $v_n$, \  and for any choice of \ $F_v{_n}$, \ there exist a choice of \ $U_{v_n}$ \ such that for any choice of \ $x_n, y_n$, \ any word along any simple path in the tree $T$ will have a length bigger than $r$.

  \bigskip

    By examining the proof of Proposition 5., we see that we have indeed proved a more general result:

  \medskip

  {\bf Proposition 6.} For any  $p\geq 4q$, property \ $T_{p,q}$ \ implies \ $TS$.

   \medskip

   {\bf Remark 7.} In the proof of Proposition 5., we had proved that $T_{4,1}\Rightarrow (TS)$.

   \medskip

   {\bf Remark 8.} In Remark 5, \ $(P')$ \ can be replaced by  \ $T_{p,q}$.

   \bigskip

    One can also generalize the property  $(P)$ by putting it in "a ternary tree form", although this generalization turns out to be equivalent to $(P)$:

    \medskip

    {\bf Property ($X$):} A finitely generated nontrivial group $K$ is said to satisfy property $(X)$, if for an infinite ternary tree with origin,  and for any \ $r>0$, there exists \ $\xi \in K$ such that for any assignment of $\xi $ and $\xi ^{-1}$ to the vertices, and any of the  assignment of elements of \ $B_r\backslash 1$ to the edges, where the elements assigned to the lower edges of the same vertex are different, any word \ $\xi ^{\pm }x_1\ldots \xi ^{\pm }x_k$ \ along any simple path has a length bigger than $r$.

    \medskip

    {\bf Remark 9.} Properties $(P)$ and $(X)$ are clearly equivalent. $(P)$ can be thought of as a "ray version" of $(X)$, but there is no difference between "ray version" and " tree version" in this case.

  \vspace{2cm}

   {\bf 3. Property $TS(\lambda )$ and product of groups.}

   \bigskip

    {\em Definition 5.} Let $\lambda \in \mathbb{R}, \lambda > 1$. A finitely generated group $\Gamma $ is said to satisfy property $TS(\lambda )$ (or belong to the class $TS(\lambda )$) if there exist $\epsilon > 0$, a finite generating set $X\subseteq \Gamma $ and $\xi \in \Gamma $ such that for any $\xi $-related set $S$, the inequality $L(S) > (\lambda + \epsilon)|S|$ is satisfied, where $L(S)$ is given in terms of the metric determined by $X$.

    \medskip

     By looking at the proof of Proposition 3, one concludes that finitely generated amenable groups never belong to the class $TS(2)$. On the other hand, if a group satisfies property $TS$ then it belongs to the class $TS(\lambda )$ for any $\lambda $. It is natural therefore to ask whether there exist groups which are not $TS$ but satisfy property $TS(\lambda )$ for some $\lambda \geq 2$.

   One can observe that products of groups are good candidates for such examples, where "mixing" the generating sets of factors can be useful to construct an example.

   \medskip

   {\bf Proposition 7.} The group $\Gamma = \mathbb{F}_2\times \mathbb{Z}$ belongs to the class $TS(\lambda )$ for any $\lambda \geq 2$.

   \medskip

   {\bf Proof.} Let $\lambda > 1$, $a, b$ be the generators of  $\mathbb{F}_2$, and $z$ be the generator of $\mathbb{Z}$. Then for any $n\in \mathbb{N}$, the set $S_n = \{ a, b, a^nz \}$ generates $\Gamma $. We will denote the Cayley metric of $\Gamma $ with respect to $S_n$ by $|.|_n$ or by $d_n(. , .)$. It is immediate that $|z^i|_n \geq n$ for all $i\in \mathbb{N}$.

 \medskip

   Let us choose $n > \lambda $ and let $\xi \in \mathbb{F}_2$ be an aperiodic word of length $4[\lambda ] + 1$. (a word $\xi \in \mathbb{F}_2$ is aperiodic if it is aperiodic as a reduced word in the alphabet $\{a, a^{-1}, b, b^{-1}\}$; in other words, if $\xi = \xi_1\omega \omega \xi _2$ for some $\xi _1, \xi _2, \omega \in \mathbb{F}_2$ and if $\omega \neq 1$ then $ |\xi | < |\xi _1| + 2|\omega | = |\xi _2|$).

   \medskip

    Assume that $S = \{(g_1, z_1), (g_1\xi , z_1), \ldots ((g_k, z_k), (g_k\xi , z_k)\}$ is a finite revised $(\xi ,0)$-related set in $\Gamma $, and let $\alpha = ((g_{n_1}, z_{n_1}), (g_{n_2}, z_{n_2}), \ldots , (g_{n_k}, z_{n_k}), (g_{n_1}, z_{n_1}))$ be the shortest closed path which passes through every element of $S$.
    Let us denote $I_1 = \{i\in \{1, \ldots ,k\} \ | \ g_{n_i} = g_{n_{i+1}}\}, I_2 = \{1, \ldots , k\}\backslash I_1$. (here, and in the sequel, $g_{n_{k+1}} = g_{n_1}$).

 \medskip

 Then $\sum _{i=1}^k d_n((g_{n_i}, z_{n_i}), (g_{n_i+1}, z_{n_i+1})) = \sum _{i\in I_1} d_n((g_{n_i}, z_{n_i}), (g_{n_i+1}, z_{n_i+1}))  +  \sum _{i\in I_2} d_n((g_{n_i}, z_{n_i}), (g_{n_i+1}, z_{n_i+1})) (\star )$

 \medskip

 Notice that $\{g_{n_i} \ | \ i\in I_2\}$ is a $\xi $-related set in $\mathbb{F}_2 = < a, b >$. By Remark 1, in the group $\mathbb{F}_2$ with a standard generating set, for any aperiodic word $\xi $ of length at least $4[\lambda ] + 1$, and for any $\xi$-related finite set $S$,  the inequality $L'(S) > \lambda |S|$ holds. Using this, from $(\star )$, we obtain that  $ \sum _{i\in I_2} d_n((g_{n_i}, z_{n_i}), (g_{n_i+1}, z_{n_i+1})) > \lambda |I_2|$. Then $\sum _{i=1}^k d_n((g_{n_i}, z_{n_i}), (g_{n_i+1}, z_{n_i+1})) > n|I_1| + \lambda |I_2| \geq \lambda (|I_1| + |I_2|) = \lambda |S|$. $\square $

  \bigskip

  {\em Definition 6.} A finitely generated group $\Gamma $ is called {\em weakly TS} if for every $\lambda > 1$ there exists a finite generating set $X\subseteq \Gamma $ and $\xi \in \Gamma $ such that for any $\xi $-related set $S$, the inequality $L(S) > \lambda |S|$ is satisfied, where $L(S)$ is defined with respect to $X$.

\medskip

 {\bf Remark 10.} Clearly, "{\em weakly TS} " is an intermediate notion between $TS$ and $TS(\lambda )$ for some $\lambda > 1$. From the proof of Proposition 9, we see that we have indeed proved more, namely, the group $\mathbb{F}_2\times \mathbb{Z}$ is weakly $TS$.

 \medskip

   We would like to raise several questions:

   \medskip

   {\bf Question 1.} If $\Gamma _1, \Gamma _2\in TS$ then is it true that $\Gamma _1\times \Gamma _2\in TS$?

   \medskip

   {\bf Remark 11.} The answer is positive if $\Gamma _1$ and $\Gamma _2$ are non-abelian free groups; the proof is similar to the proof of Proposition 2 in [Ak2].

    \medskip

  {\bf Conjecture 1.}  R.Thompson's group $F$([CFP]) belongs to the class $TS(2)$.

   \medskip

  {\bf Question 2.} Is there a group which is $TS(\lambda )$ for some $\lambda \geq 2$ but is not weakly $TS$?

    \medskip

   {\bf Remark 12.} It is proved in [Gu] that $F$ does not satisfy  the inequality of the condition $TS(2.5)$ with respect to standard generating set.
This disproves my original conjecture that $F$ might be even a $TS = TS(\infty )$ group. Conjecture 1. seems to be the right one.

   \bigskip

   It was suggested by W.Thurston that some modification of the notion $TS$ may be equivalent to non-amenability. The notions of {\em weakly TS} and $TS(\lambda )$ groups could be the beginning of the right approach to this question.

   But $TS(\lambda )$ groups turn out very useful in the areas of application, especially in the theory of non-amenable groups.

   Below we introduce properties $P(r), P(r)', P_n(r)'$ where $r\in \mathbb{R}_+, n\in \mathbb{N}$. Let $K$ be a finitely generated nontrivial group with a fixed left-invariant metric.

   \medskip

    $(P(r))$:  $K$ is said to satisfy the property $P(r)$ if there exists \ $\xi \in K$ \ such that $$|\xi ^{\epsilon _1 }x_1\xi ^{\epsilon _2 }x_2\ldots \xi ^{\epsilon _k }x_k| > r $$ for any positive $k$ and \ $x_1,\ldots ,x_k \in K, \ x_i \neq 1, \ |x_i| \leq r, \ \epsilon _i\in \{+1, -1\}$

\medskip

 $(P(r)')$:  $K$  is said to satisfy a property $(P(r)')$ , if there exists \ $\xi \in K$ \ such that for every aperiodic sequence \ $x_1, x_2, \ldots x_k$ of elements in  $K$, the inequality \ $|\xi ^{\epsilon _1}x_1\ldots \xi ^{\epsilon _k}x_k| > r$ \ satisfied for every \ $k\in \mathbb{Z}_+, \ x_i \neq 1, \ |x_i| \leq r, \ \epsilon _i \in \{+1, -1\}$.

  \medskip

  $(P_n(r)')$. Let $n$ be a positive integer.  $K$  is said to satisfy a property $(P_n(r)')$ if there exists \ $\xi \in K$ \ such that for every $n$-aperiodic sequence \ $x_1, x_2, \ldots x_k$ of elements in \ $K$, the inequality \ $|\xi ^{\epsilon _1}x_1\ldots \xi ^{\epsilon _k}x_k| > r$ \ satisfied for every \ $k\in \mathbb{Z}_+, \ x_i \neq 1, \ |x_i| \leq r, \ \epsilon _i \in \{+1, -1\}$.

  \bigskip

   The proof of the following proposition is identical to the proof Proposition 5. except wherever we use property $P$ or $P'_{10}$, we need to use property $P(r)$ or $P'_{10}(r)$.

   \medskip

   {\bf Proposition 8.} For a finitely generated group $\Gamma $ with a fixed left-invariant metric, \ $P(r)\Rightarrow TS(r/12)$ and $P_{10}(r)'\Rightarrow TS(r/96)$.

   \medskip

   {\bf Proposition 9.} A finitely generated group satisfying the condition $P_{10}'(192)$ is non-amenable.

   \medskip

   {\bf Proof.} By Proposition 8. we obtain that $P_{10}'(192)\Rightarrow TS(2)$. But $TS(2)$ implies non-amenability. $\square $

 \vspace{1.3cm}

 {\bf 5. Applications.}

  \bigskip

  Below, we will discuss some applications of our criteria.

  \medskip

  First, let us mention an example in which the direct application of the criteria fails immediately:

   \medskip

  {\em Free object of the variety $[X^p,Y^p] = 1$ where $p$ is any positive integer:} \ Let $G_{n,p}$ be a free object on $n\geq 1$ generators. If $n=1$ then $G_{n,p}$ is isomorphic to $\mathbb{Z}$; for $n\geq 2$, assume $x_i = u_i^p, \ i\in \{1, 2, \ldots , 2k\}$ such that the sequences $(x_1, x_3, \ldots x_{2k-1})$ and $(x_2, x_4, \ldots x_{2k})$ satisfy the following conditions: for every $i\in \{1, 2, \ldots , 2k\}, \  Card \{ j \ | \ x_{2j-1} = u_i^p\}  =  Card \{ j \ | \ x_{2j-1} = u_i^{-p}\}$ and  $Card \{ j \ | \ x_{2j} = u_i^p\}  =  Card \{ j \ | \ x_{2j} = u_i^{-p}\}$.

   Then for every $\xi \in G_{n,p}$, $$\xi x_1\xi ^{-1} x_2\xi x_3\xi ^{-1}x_4\ldots \xi x_{2k-1}\xi ^{-1}x_{2k} = (\xi u_1\xi ^{-1})^p u_2^p (\xi u_3\xi ^{-1})^pu_4^p\ldots  (\xi u_{2k-1}\xi ^{-1})^pu_{2k}^p = 1$$ Since, for any $m\in \mathbb{N}$, it is possible to choose $k\in \mathbb{N}$ and $m$-aperiodic sequence $(x_1, \ldots , x_{2k})$ satisfying the above conditions, the group $G_{n,p}$ does not satisfy property $(P_m')$ or $(T_{p,q})$.

   \medskip

   {\bf Application to group varieties.} Since verifying the condition $(P_{10}(192)')$ immediately yields non-amenability, this allows to prove non-amenability for free objects (relatively free groups) of some group varieties. The most  interesting example is, of course, the free Burnside group $\mathbb{B}(m,n)$ for $m\geq 2$ and $n\geq 665$ with $n$ being odd.

    \medskip

 The condition $(P_{10}(192)')$ is indeed satisfied for free Burnside groups of sufficiently large odd exponent. First, we need the following two lemmas:

 \medskip

 {\bf Lemma 4.} There exists a reduced word $\xi $ in the alphabet $S = \{a, b, a^{-1}, b^{-1}\}$ such that

    (i) $\xi $ is 3-aperiodic.

  (ii) $|\xi | \geq 10000$

  (iii) $\xi $ satisfies condition $C'(1/5)$.

  (iv) $\xi \equiv \alpha \xi ^{(0)}\beta $ where $|\alpha | = |\beta | = 400$ and the set $\{\alpha , \beta \}$ (i.e. the symmetrization of this set) satisfies condition $C'(1/3)$.

  \medskip

  Proof. For the definition of the small cancelation conditions we refer the reader to [G-H], page 228.
Let $x = b, y = aba, z = aabaa$. Let $\delta $ be an aperiodic sequence in $x, y, z$ such that the length of $\delta $ in the alphabet of $S$ is equal to $N\in (10000, 10006)$. Then $\delta $ is 3-aperiodic in $S$.

  Let $\delta = c_1c_2\ldots c_N$ where $c_i\in S, 1\leq i\leq N$. Denote $B = \{j \ | \ 1\leq j\leq N, c_j = b\}$. Notice that $B$ intersects every subinterval $[u, u+4]\subseteq [1, N], u\in \mathbb{N}$ of length 4. Let $B = \{b_1, b_2, \ldots , b_m\}$ where the elements are listed in the increasing order. We choose a subset $D = \{d_1, d_2, \ldots , d_t\} \subset B$ with elements again listed in the increasing order (i.e. $d_1 < d_2 < \ldots < d_t$) such that the following conditions are satisfied:

 (1) $d_{p+1}-d_p \neq d_{q+1}- d_q$ for all distinct $p, q\in \{1, 2, \ldots , t-1\}$

 (2) $D\cap [u, u+60] \neq \emptyset $ for all $u\in (0, 340)\cup (N-400, N-60)$.

 (3) $D\cap [u, u+1000]\neq \emptyset $ for all $u\in (0, N-1000)$.

 (4) $|D\cap (B\cap [u, u+100])| < \frac {1}{3}|B\cap [u, u+100]|$ for all $u\in (0, N-100)$.

  Then by replacing all letters $b$ by $b^{-1}$ in the index set $D$ (i.e. for all $i\in D$, $c_i = b$ is replaced by $b^{-1}$), we obtain from $\delta $ a new word $\xi $ satisfying all the conditions (i)-(iv). $\square $

  \bigskip

  {\bf Lemma 5.} If the sequence $x_1, x_2, \ldots , x_k$ of non-trivial reduced words of length at most 192 in the alphabet of standard generating set of the free group $\mathbb{F}_n, n\geq 2$ is 10-aperiodic then we can choose $\xi \in \mathbb{F}_n$ such that a word  $ \xi ^{\epsilon _1}x_1\ldots \xi ^{\epsilon _k}x_k, \ \epsilon _1, \ldots \epsilon _k\in \{-1, 1\}$, in its reduced form is $500$-aperiodic

 \medskip

 Proof. Let $\xi \in \mathbb{F}_n$ be a reduced word satisfying the following conditions:

  (i) $\xi $ is 3-aperiodic.

  (ii) $|\xi | \geq 10000$

  (iii) $\xi $ satisfies condition $C'(1/5)$.

  (iv) $\xi \equiv \alpha \xi ^{(0)}\beta $ where $|\alpha | = |\beta | = 400$ and the set $\{\alpha , \beta \}$ satisfies condition $C'(1/3)$.

  Because of conditions (ii) and (iv), for all $k\in \mathbb{N}$, and for all $\epsilon _i\in \{-1, 1\}, i\in \{1, 2, \ldots , k\}$,    we have $ \xi ^{\epsilon _1}x_1\ldots \xi ^{\epsilon _k}x_k = \xi _1u_1\xi _2u_2\ldots \xi _ku_k$ where $\xi _1u_1\xi _2u_2\ldots \xi _ku_k$ is a reduced word, moreover, $\forall  i\in \{1, 2, \ldots , k\}, \ \xi _i$ is a subword of $\xi $ or $\xi ^{-1}$, and $|\xi _i| \geq |\xi | - 400, |u_i| \leq 400$.

  \medskip

  Assume that the reduced word   $\xi _1u_1\xi _2u_2\ldots \xi _ku_k$ is not 500-aperiodic and let $A$ be a subword of it of order at least 500.

  We will consider two cases:

  \medskip

  {\em Case 1.} $|A| \leq \frac {1}{5}|\xi |$.

   In this case, since $|\xi _i| > \frac {4}{5}|\xi |$ and $|u_i|\leq 400$, for all $i\in \{1, 2, \ldots , k\}$, we obtain that $AAAA$ is a subword of $\xi _i$ for some $i\in \{1, 2, \ldots , k\}$ therefore $AAAA$ is a subword of $\xi $ or $\xi ^{-1}$ which contradicts condition (i).

   \medskip

   {\em Case 2.}  $|A| \geq \frac {1}{5}|\xi |$

  Notice that for all $i\in \{1, 2, \ldots , |A|\}$, a shift $R_i(A)$ of $A$ is also a period of order at least $500-1 = 499$. Then there exists $j\in \{1, 2, \ldots , k\}$ (even more precisely - because of inequality  $|A| \geq \frac {1}{5}|\xi |$ - there exists $j\in \{1, 2, \ldots , k-100\}$) and a subword $B$ of period at least 499 such that $\xi _1u_1\xi _2u_2\ldots \xi _ku_k \equiv \xi _1u_1\ldots \xi _ju_j\xi _{j+1}'B^{499}X$ where $\xi _{j+1}'$ is a subword of $\xi _{j+1}$ and $|\xi _{j+1}|'\in (\frac {1}{2}|\xi | - 500,  \frac {1}{2}|\xi | + 500)$.

  Then either condition (iii) is violated or the sequence $x_{j+1}, x_{j+2}, \ldots , x_k$ is not 499-aperiodic in which case it is not 10-aperiodic which contradicts the assumption of the lemma. $\square $

 \bigskip

 Now, for the group $\mathbb{B}(m,n)$ for $m\geq 2$, fixing the standard generating set in this group, if $x_1, x_2, \ldots , x_k$ is a 10-aperiodic sequence where every $x_i, 1\leq i\leq k$ is represented by the shortest word in the standard generating set, and if $x_i \in B_{192}\backslash \{1\}$, then, by Lemma 6, we can choose a word $\xi $ such that a word $ \xi ^{\epsilon _1}x_1\ldots \xi ^{\epsilon _k}x_k$ is $500$-aperiodic. From the proof of Theorem 19.1. in [Ol3], we see that  if the exponent of the Burnside group is sufficiently big and odd, then 500-aperiodic words are all distinct thus the inequality  $|\xi ^{\epsilon _1}x_1\ldots \xi ^{\epsilon _k}x_k| > 192$ is satisfied.

    \bigskip
    
    {\bf Remark 13.} Although very difficult to prove, but it is one of the  basic facts of the theory of free Burnside groups that in $\mathbb{B}(m, n)$, where $m\geq 2$ and $n\in \mathbb{Z}_{odd}$ is sufficiently big, aperiodic words (or $k$-aperiodic words) are distinct. It is this fact which immediately implies that these groups are infinite, and even more, they have exponenetial growth. Notice that we are using just this very basic fact from the theory of Burnside groups to establish their non-amenability as well.   
    
    \medskip
    
    {\bf Remark 14.} In [Ol2], the author proves that, for all odd $n > 10^{10}$, $[n/3]$-aperodic words in $\mathbb{B}(2,n)$ are non-trivial.
Notice that this (more advanced) fact also suffices for our purposes, in fact, it is much more than what we need.

  \medskip
  
    {\bf Remark 15.} It is true that even for exponent $n = 665$ the condition $(P)'$ is satisfied (private communication with S.Adian). Although we only have proved that $(P_{10})'\Rightarrow TS\Rightarrow ${\em non-amenability} the proof can be modified to obtain $(P')\Rightarrow TS\Rightarrow ${\em non-amenability}, which then implies non-amenability even for the exponent 665. Finally, notice that non-amenability of relatively free group of the variety $X^p = 1$ implies non-amenability of relatively free groups of the variety $U(X^p, Y^p) = 1$ where $U$ is any non-trivial word in the free group on 2 generators.  

    \vspace{3cm}

    { REFERENCES:}

    \bigskip

   [Ad1] S.Adian, \ Random walks on free periodic groups. (Russian) Izv. Akad. Nauk SSSR Ser. Mat. 46 (1982), no. 6, 1139--1149, 1343.

   [Ad2] S.Adian, \ The Burnside problem and identities in groups. \ Ergebnisse der Mathematik und Ihrer Grenzgebeite.

    [Ak1] Akhmedov, A. Ph.D Thesis. \ Yale University, 2004.

    [Ak2] Akhmedov, A. Traveling Salesman Problem in Groups.  {\em Contemp. Math.}, vol. 372, Amer. Math. Soc., Providence, RI, 2005.

    [Ak3] Akhmedov, A. Perturbations of Wreath Product and Quasi-Isometric Rigidity. {\em IMRN 2008.}

    [CFP] Cannon,J.W., Floyd,W.J., Parry,W.R. \ Introductory notes on Richard Thompson's groups. \ {\em Enseign.  Math.}  (2)  {\bf 42}  (1996), \ no 3-4.

    [G-H] E.Ghys, P, de la Harpe, \ Sur les Groupes Hyeprboliques d'apres Mikhael Gromov.  \ Birkhauser 1990.

    [Gr] Grigorchuk,R. \ Symmetrical random walks on discrete groups. {\em Multicomponent random systems}, pp. 285--325, \ 1980.
    
    [Gu] Guba, V. Traveller Salesman Property and R.Thompson's Group F. 137Ð142, {\em Contemp. Math.}, 394, Amer.  Math. Soc., Providence, RI, 2006.

    [Ke] Kesten, H. Full Banach mean values on countable groups. Math. Scand. 7 (1959), 146{156.
    
   [Kr] Kruskal, J.  \ On the shortest spanning subtree of a graph and the traveling salesman problem. \ {\em Proc. Amer. Math. Soc.} {\bf 7}, 1956.

   [Ol1] Olshanski, A. \ On the question of the existence of an invariant mean on a group. (Russian) \ {\em Uspekhi Mat. Nauk.} {\bf 35.} (1980), no. 4(214), 199-200;

   [Ol2] Olshanski, A. \ On Novikov-Adian Theorem. \ {\em Math. Sbornik.} {\bf 118}, (1982),no, 203-235
   
   [Ol3] Olshanski, A. \ The geometry of defining relations in groups, Moscow: Nauka, 1989.
   
   [OS] Olshanski, A., Sapir, M. \ Non-amenable finitely presented torsion-by-cyclic groups. Publ. Math. Inst. Hautes ƒtudes Sci. No. {\bf 96} (2002), 43--169 (2003)

   \end{document}